\documentclass[a4paper,10pt,leqno]{amsart}

        \title[Coefficients for the Farrell-Jones Conjecture]
              {Coefficients for\\ the Farrell-Jones Conjecture}
       \author{Arthur Bartels}
       \author{Holger Reich}
      \address{Westf\"alische Wilhelms-Universit\"at M\"unster\\
               Mathematisches Institut\\
               Einsteinstr.~62,
               D-48149 M\"unster, Germany}
        \email{bartelsa@math.uni-muenster.de}
      \urladdr{http://www.math.uni-muenster.de/u/bartelsa/bartels}
        \email{reichh@math.uni-muenster.de}
      \urladdr{http://www.math.uni-muenster.de/u/reichh}
         \date{\today}
    \subjclass[2000]{}

\usepackage[
]{hyperref}
\usepackage{enumerate,amssymb}
\usepackage[arrow,curve,matrix,tips,2cell]{xy}
  \SelectTips{eu}{10} \UseTips
  \UseAllTwocells

\DeclareMathAlphabet{\matheurm}{U}{eur}{m}{n}

\newcommand{\hf}{\text{hf}}
\newcommand{\fr}{\text{fr}}

\newcommand{\pt}{\text{pt}}

\newcommand{\Ab}{\matheurm{Ab}}
\newcommand{\AddCat}{\matheurm{Add~Cat}}
\newcommand{\AddCatInv}{\matheurm{Add~Cat~Inv}}

\newcommand{\Or}{\matheurm{Or}}

\newcommand{\Sp}{\matheurm{Sp}}

\DeclareMathOperator{\Aut}{Aut}

\DeclareMathOperator{\ch}{ch}

\DeclareMathOperator{\card}{card}

\DeclareMathOperator{\cone}{cone}
\DeclareMathOperator{\cyl}{cyl}

\DeclareMathOperator{\id}{id}
\DeclareMathOperator{\inc}{inc}
\DeclareMathOperator{\ind}{ind}

\DeclareMathOperator{\map}{map}
\DeclareMathOperator{\mor}{mor}

\DeclareMathOperator{\res}{res}

\DeclareMathOperator{\supp}{supp}
\DeclareMathOperator{\Sw}{Sw}

\DeclareMathOperator*{\sma}{\wedge}  

\newcommand{\Fin}{{\mathcal{F}\text{in}}}

\newcommand{\VCyc}{{\mathcal{VC}\text{yc}}}

  \newcommand{\IK}{\mathbb{K}}
  \newcommand{\IL}{\mathbb{L}}
  
  \newcommand{\IN}{\mathbb{N}}

  \newcommand{\IR}{\mathbb{R}}
  
  \newcommand{\IT}{\mathbb{T}}

  \newcommand{\IZ}{\mathbb{Z}}

  \newcommand{\cala}{\mathcal{A}}
  \newcommand{\calb}{\mathcal{B}}
  \newcommand{\calc}{\mathcal{C}}
  \newcommand{\cald}{\mathcal{D}}
  \newcommand{\cale}{\mathcal{E}}
  \newcommand{\calf}{\mathcal{F}}
  \newcommand{\calg}{\mathcal{G}}

  \newcommand{\bfE}{{\mathbf E}}
  
  \newcommand{\bfF}{{\mathbf F}}

  \newcommand{\bfK}{{\mathbf K}}
  \newcommand{\bfL}{{\mathbf L}}

\newcommand{\x}{{\times}}
\newcommand{\ox}{{\otimes}}

\theoremstyle{plain}      \newtheorem{theorem}{Theorem}[section]
                          
                          \newtheorem{lemma}[theorem]{Lemma}
                          \newtheorem{corollary}[theorem]{Corollary}
                          \newtheorem{proposition}[theorem]{Proposition}
                          
                          \newtheorem{conjecture}[theorem]{Conjecture}

\theoremstyle{definition} \newtheorem{definition}[theorem]{Definition}
                          \newtheorem{example}[theorem]{Example}

                          \newtheorem{remark}[theorem]{Remark}

\theoremstyle{remark}

\makeatletter\let\c@equation=\c@theorem\makeatother

\begin{document}

\begin{abstract}
We introduce the Farrell-Jones Conjecture with coefficients in
an additive category with $G$-action.
This is a variant of the Farrell-Jones Conjecture
about the algebraic $K$- or $L$-Theory of a group ring $RG$.
It allows to treat twisted group rings and crossed product rings.
The conjecture with coefficients
is stronger than the original conjecture but it
has better inheritance properties.
Since known proofs using controlled algebra
carry over to the set-up with coefficients we obtain new results
about the original Farrell-Jones Conjecture.
The conjecture with coefficients implies the
fibered version of the Farrell-Jones
Conjecture.
\end{abstract}

\maketitle


\section{Introduction}
\label{sec:introduction}

The Farrell-Jones Conjecture predicts that
the algebraic $K$- or $L$-theory of a group ring $RG$ can be
described in terms of the $K$- respectively $L$-theory of group rings $RH$, where $H$ ranges over
the  family of virtually cyclic subgroups of $G$, compare \cite{FJ-iso-conj}.
More formally the conjecture says that the assembly map
\[
H_*^G ( E_{\calf} G ; \bfK_R ) \to K_* ( RG ),
\]
which assembles $K_*(RG)$ from the pieces $K_* (RH)$, $H \in \calf$,
is an isomorphism if $\calf$ is the family of virtually cyclic subgroups of $G$.
Here $H_{\ast}^G ( - ; \bfK_R )$ is a $G$-equivariant homology theory and $E_{\calf} G$ denotes
the classifying space for the family of subgroups $\calf$. A group is called virtually cyclic if
it contains a cyclic subgroup of finite index. For more explanations see \cite{Lueck-Reich-survey}.
There is a similar formulation for $L$-theory.

The goal of this paper is to define, and prove in many cases, a
``Farrell-Jones Conjecture with coefficients in an additive category with $G$-action''.
A precise formulation is given in Section~\ref{sec:FJ-with-coeff} and makes essential use of Definition~\ref{def:cala-ast-T}.
The conjecture with coefficients  is a generalization of the original conjecture.
The special case
where the additive category with $G$-action is the category of finitely generated free $R$-modules
equipped with the trivial $G$-action corresponds to
the usual Farrell-Jones Conjecture.

The reason for considering this sort of generalization of the Farrell-Jones Conjecture is twofold.
First, the conjecture with coefficients
has better inheritance properties.
\begin{theorem}
The Farrell-Jones Conjecture with coefficients in an additive category with $G$-action passes to arbitrary subgroups and more
generally it ``pulls back'' under arbitrary group homomorphisms, where of course the family needs to be pulled back as well.
But even more is true: the injectivity and surjectivity part of the conjecture have these inheritance properties independently.
\end{theorem}
For a precise statement see Conjecture~\ref{con:FJC-with-coefficients-alg-K-theory},
Corollary~\ref{cor:Inheritance-for-group-homoeomorphisms} and Theorem~\ref{thm:FJ-w-C-passes-to-subgroups-alg-K}.

Second, we obtain as a special case of the conjecture with coefficients the
correct conjecture for the $K$-theory of twisted group rings and more generally
crossed product rings. For example let
$G$ operate through ring homomorphisms on the ring $R$,
i.e.\ we have a group homomorphism $\alpha \colon G \to \Aut(R)$, and let $R_{\alpha} G$ denote the twisted group ring.
Then  of course one expects
that $K_*( R_{\alpha} G)$ can be assembled from the $K_* ( R_{\alpha} H )$,
where $H$ ranges over the virtually cyclic subgroups of $G$.
Our conjecture makes this precise, compare Conjecture~\ref{con:FJ-for-crossed-product}.

Recall that crossed product rings
play an important role in Moody's Induction Theorem, see \cite{Moody(1987)} and Chapter~8 in \cite{Passman(1989)}.

Of course all this is only useful if there are techniques which prove the more general conjecture.
Many results on the Farrell-Jones Conjecture (without coefficients)
use the concept of controlled algebra and the description of the assembly map as a ``forget-control map''.
In Section~\ref{sec:controlled-algebra} we extend these concepts to the case with coefficients and
formulate a ``forget-control version'' of the Farrell-Jones Conjecture with Coefficients.
By simply inspecting existing proofs in the literature, see Section~\ref{sec:Applications},
one obtains results about the Farrell-Jones Conjecture with Coefficients. Combined with the inheritance properties
we obtain the following new results about the original Farrell-Jones Conjecture.

\begin{corollary}
Let $G$ be the fundamental group of a
closed Riemannian manifold
of strictly negative sectional curvature
and let $\Gamma$ be a subgroup of $G$.
Then for every ring $R$ the assembly map
\[
H^\Gamma_* (E_\VCyc \Gamma; \bfK_R) \to K_* (R\Gamma)
\]
is an isomorphism.
\end{corollary}
This is an extension of the main result from \cite{BR-JAMS}.
It follows from Theorem~\ref{thm:FJ-w-C-passes-to-subgroups-alg-K},
Remark~\ref{rem:with-implies-without} and Theorem~\ref{thm:BR-with-coefficients}.

\begin{corollary}
Let $\Gamma$ be a subgroup of a hyperbolic group
in the sense of Gromov
and let $R$ be a ring.
Then the assembly map
\[
H^\Gamma_*(E_\Fin \Gamma; \bfK_R) \to K_*(R\Gamma)
\]
is split injective.
\end{corollary}

This generalizes \cite{Rosenthal-splitting-K-theory} and uses the fact
proven in \cite{Rosenthal-Schuetz-hyperbolic} that
hyperbolic groups in the sense of Gromov
satisfy the assumptions of Theorem~\ref{thm:david-with-coefficients}.

\begin{corollary}
Let $G$ be group of finite asymptotic dimension
that admits a finite model for the classifying space $BG$.
Let $\Gamma$ be a subgroup of $G$.
Then for every ring $R$ the assembly map
\[
H^{\Gamma}_*(E\Gamma ;\bfK_R) \to K_*(R \Gamma)
\]
is split injective.
\end{corollary}

This is a generalization of \cite{Bartels-finite-asymp-dim} and follows because of the inheritance properties immediately from
Theorem~\ref{thm:bartels-fin-asy-with-coeff}.

In \cite{FJ-iso-conj} Farrell and Jones develop the Fibred Isomorphism Conjecture,
a different generalization of the Farrell-Jones Conjecture,
which also has better inheritance properties.
The fibred version is however not so well adapted to proofs which use controlled algebra as opposed to
controlled topology.
The precise relationship between the Fibered Farrell-Jones Conjecture and the Farrell-Jones Conjecture with Coefficients
is discussed in Remark~\ref{rem:coefficients-implies-fibered}.

In the context of  topological $K$-theory of $C^{\ast}$-algebras there is an analog to the Farrell-Jones
Conjecture with
Coefficients~\ref{con:FJC-with-coefficients-alg-K-theory},
the Baum-Connes Conjecture with Coefficients. Of course the
development of a Farrell-Jones Conjecture with Coefficients was motivated by this analogy.
The reader should be warned that the
Baum-Connes Conjecture with Coefficients
is know to be wrong
\cite{Higson-Lafforgue-Skandalis-counterexamples}.
At present it is not clear whether the
Farrell-Jones Conjecture with
Coefficients~\ref{con:FJC-with-coefficients-alg-K-theory}
fails for the groups considered in
\cite{Higson-Lafforgue-Skandalis-counterexamples}.

We would like to thank David Rosenthal for helpful discussions on his paper \cite{Rosenthal-splitting-K-theory}.
\medskip

The paper is organized as follows\,:

\begin{tabular}{ll}
\ref{sec:introduction}. & Introduction
\\
\ref{sec:the-category}. & The category $\cala \ast_G T$
\\
\ref{sec:FJ-with-coeff}. & The Farrell-Jones Conjecture with Coefficients
\\
\ref{sec:inheritance-properties}. &  Inheritance properties
\\
\ref{sec:L-theory}. & $L$-theory
\\
\ref{sec:crossed-products}. & Crossed Products
\\
\ref{sec:controlled-algebra}. & Controlled Algebra
\\
\ref{sec:Applications}. &  Applications
\\
\ref{sec:swan-group-action}. & The Swan group action
\\
 & References
\end{tabular}

\section{The category $\cala \ast_G T$}
\label{sec:the-category}

In the following we will consider additive categories $\cala$ with a
right $G$-action, i.e.\
to every group element $g$ we assign an additive covariant functor
$g^{\ast} \colon \cala \to \cala$,
such that $1^{\ast} = \id$ and composition
of functors (denoted $\circ$) relates to multiplication in the group
via $g^{\ast} \circ h^{\ast} = (hg)^{\ast}$.

\begin{definition}
\label{def:cala-ast-T}
Let $\cala$ be an additive category with a right $G$-action
and let $T$ be a left $G$-set.
We define a new additive category denoted $\cala \ast_G T$ as follows.
An object $A$ in $\cala \ast_G T$ is a family
\[
A = ( A_t) _{t \in T}
\]
of objects in $\cala$ where we require
that $\{ t \in T \; | \; A_t \neq 0 \}$ is a finite set.
A morphism $\phi \colon A \to B$ is a collection of morphisms
\[
\phi = (\phi_{g,t} )_{(g , t) \in G \times T},
\]
where
\[
\phi_{g,t} \colon A_t \to g^{\ast} \left( B_{gt} \right)
\]
is a morphism in $\cala$. We require that the set of
pairs $(g,t) \in G \x T$  with $\phi_{g,t} \neq 0$ is finite.
Addition of morphisms is defined componentwise.
Composition of morphisms is defined as follows.
Let $\phi = ( \phi_{g,t} ) \colon A \to B$ and
$\phi^{\prime} = ( \phi^{\prime}_{g,t} ) \colon B \to C$
be given then the composition
$\psi = \phi^{\prime} \circ \phi \colon A \to C$ has components
\begin{equation}
\label{eq:composition-in-cala-ast}
\psi_{g,t} = \sum_{h, k \in G \quad g = kh}
h^{\ast} ( \phi^{\prime}_{k , ht} ) \circ \phi_{h,t}.
\end{equation}
\end{definition}

The reader could now pass immediately to Section~\ref{sec:FJ-with-coeff} in order to see how
the Farrell-Jones Conjecture with Coefficients is formulated.

\begin{remark}[Naturality of $\cala \ast_G T$]
\label{rem:cala-ast-T-is-natural}
The construction $\cala \ast_G T$ is natural in $\cala$,
i.e.\ if $F \colon \cala \to \cala'$ is an additive functor which is
equivariant with respect to the $G$-action then
\[
(F \ast_G T(A))_t = F(A_t) \quad \mbox{and}
      \quad (F \ast_G T(\phi))_{g,t} = F(\phi_{g,t})
\]
defines a functor
$F \ast_G T \colon \cala \ast_G T \to \cala^{\prime} \ast_G T$.
If the functor $F$ is an equivalence of categories, then $F \ast_G T$ is
an equivalence of categories.

If $f \colon T \to T^\prime$ is a $G$-equivariant map then
\[
(\cala \ast_G f(A))_{t^\prime} = \oplus_{t \in f^{-1}(t^\prime)} A_t
      \quad \mbox{and} \quad (\cala \ast_G f(\phi))_{g,t^\prime} =
             \oplus_{t \in f^{-1}(t^\prime)} \phi_{g,t}
\]
defines (almost) a functor
$\cala \ast_G f \colon \cala \ast_G T \to \cala \ast_G T^\prime$.
The minor problem, that this definition involves the choice of
a direct sum, can be resolved by redefining an object
in $\cala \ast_G T$ to be an object $A$ as before together with
a choice of direct sum $\oplus_{s \in S} A_s$ for every subset $S$ of
$T$. We prefer to ignore this problem.
\end{remark}

\begin{example}[Trivial action]
\label{ex:cala-ast-T-and-Davis-Lueck}
A ring $R$ can be considered as a category with one object.
Let $R_{\oplus}$ denote the additive category obtained from $R$ by
formally adding sums, compare
\cite[Exercise~5 on p.194]{MacLane-cat-for-working}.
This is a small model for
the category of finitely generated free $R$-modules.
If we equip $R_\oplus$ with the trivial right $G$-action
then there is a canonical identification
\[
R_{\oplus} \ast_G T = R \calg^{G} (T)_{\oplus},
\]
where the right hand side is constructed in
\cite[Section 2]{Davis-Lueck-assembly}
and plays an important role in the construction
of the assembly map by Davis and L\"uck.
This identification is natural with respect to maps of $G$-sets
$T \to T'$. Let us specialize to the case where $T=G/H$. Then the inclusion of the
full subcategory consisting of objects $A=(A_{gH})$, with $A_{gH} = 0$ for $gH \neq eH$, induces an equivalence
\begin{eqnarray} \label{eq-FRH-R-stern-pt}
\calf^f( RH ) \simeq R_{\oplus} \ast_G G/H .
\end{eqnarray}
Here $\calf^f (RH)$ denotes the category of finitely generated free $RH$-modules.
\end{example}

\begin{example}[Twisted group rings] \label{ex-twisted-group-rings}
Suppose the group $G$ acts via ring homomorphisms on $R$, i.e.\ we are given a group homomorphism
$\alpha \colon G \to \Aut(R)$. Then the twisted group ring $R_{\alpha}G$ is defined to be $RG$ as an abelian group with
the twisted multiplication determined by $g r = \alpha(r) g$ for $r \in R$ and $g \in G$.
There is a right $G$-operation defined on the category of $R$-modules, where $g^{\ast}M = \res_{\alpha(g)} M$, i.e.\
$g^{\ast}M$ has the same underlying abelian group but the $R$-module structure is twisted by $\alpha$, i.e.\
$r \cdot_{\res_{\alpha(g)} M} m = \alpha(g)(r) m$ for $r \in R$ and $m \in M$.
Let $\calf^f(R)$ denote a small model for the category of finitely generated free right $R$-modules.
One can arrange that $\calf^f(R)$ is closed under
the $G$-operation. We show in Section~6 that there is an equivalence of categories
\[
\calf^f(R) \ast_G \pt \simeq \calf^f ( R_{\alpha} G )
\]
and that more generally
\[
\calf^f(R) \ast_G G/H \simeq \calf^f(R_{\alpha|_H} H ).
\]
\end{example}

\begin{example}[Group extensions] \label{ex-crossed-product}
Suppose $K$ is a normal subgroup of $\Gamma$ and let $p \colon \Gamma \to \Gamma/K=G$ denote the quotient homomorphism.
If the group extension
$1 \to K \to \Gamma \to G \to 1$ splits we can choose a group homomorphism $s \colon G \to \Gamma$ such that $p \circ s  = \id$.
If we define $\alpha(\gamma) \colon RK \to RK$
as conjugation with $s(\gamma)$ we see that the group ring $R\Gamma$ can be written as a twisted group
ring $RK_{\alpha} G$, compare Example~\ref{ex-twisted-group-rings}.
If however the extension is non-split and $s$ is only a set-theoretical section (with $s(1)=1$) then
the group ring $R\Gamma$ is a crossed product ring
\[
R\Gamma = RK_{\alpha , \tau} G,
\]
and no longer a twisted group ring,
compare~\cite{Passman(1989)} and Section~\ref{sec-crossed} below.
In particular
$\gamma \mapsto \alpha ( \gamma )$ no longer defines an action of $G=\Gamma/K$ on $RK$ and
correction terms, expressible in terms of
the cocycle $\tau(\gamma , \gamma^{\prime}) = s( \gamma ) s( \gamma^{\prime}) s( \gamma \gamma^{\prime} )^{-1}$ will
appear.

The language developed above absorbs these extra difficulties. We will see that we can work with an honest action
if we are working in the context of actions on additive categories.
By \eqref{eq-FRH-R-stern-pt} we have
\[
\calf^f(RK) \xrightarrow{\simeq} R_{\oplus} \ast_\Gamma \Gamma/K.
\]
The category on the right should be thought of as a ``fat'' version of the
category $\calf^f( RK )$ of finitely generated free $RK$-modules,
which has the advantage that it carries an honest naturally defined right $\Gamma/K$-action.
Now applying $- \ast_{\Gamma/K} \pt$ should be thought of as forming a ``fattened'' twisted group ring,
compare Example~\ref{ex-twisted-group-rings}.
In fact
we have an additive equivalence
\begin{eqnarray*}
(R_{\oplus} \ast_\Gamma \Gamma/K ) \ast_{\Gamma/K} \pt  & \xrightarrow{\cong} &
R_{\oplus} \ast_{\Gamma \times \Gamma/K} ( \Gamma/K \times \pt ) \\
 & \xrightarrow{\simeq} & R_{\oplus} \ast_\Gamma \pt = R\Gamma_{\oplus}
\end{eqnarray*}
by an application of Proposition~\ref{prop:main-technical-ast-properties}\ref{ast-first}
respectively \ref{ast-second} below.
\end{example}

Let $K$ and $G$ be groups.
If $\cala$ is an additive category with right $K$-action and
$S$ is a $K$-$G$ biset.
Then $\cala \ast_K S$ can be equipped with a right $G$-action as follows.
If $A = (A_s )_{s \in S}$ is an object
and $\phi = ( \phi_{k,s} )_{(k,s) \in K \times S}$ is a morphism in
$\cala \ast_K S$ then for $g \in G$
we define $g^{\ast} A$ and $g^{\ast} \phi $ by
\[
(g^{\ast}  A)_s = A_{sg^{-1}}  \quad \mbox{ and }
\quad (g^{\ast}  \phi )_{k,s} = \phi_{k , sg^{-1}}.
\]

\begin{proposition} \label{prop:main-technical-ast-properties}
\begin{enumerate}
\item \label{ast-first}
Let $K$ and $G$ be groups.
Suppose $\cala$ is an additive category with right $K$-action,
let $S$ be a $K$-$G$ biset and let $T$ be a left $G$-set.
Then there is an additive isomorphism of additive categories
\[
( \cala \ast_K S ) \ast_G T \xrightarrow{\cong}
             \cala \ast_{K \times G} ( S \times T ).
\]
Here in order to form the category on the right hand side we let
$K \times G$ act from the right on $\cala$ via the projection to $K$,
and from the left on the set $S \times T$ by
$(k,g)(s,t) = (ksg^{-1} , gt)$ for
$(k,g) \in K \times G$ and $(s,t) \in S \times T$.
\item \label{ast-second}
Let $N$ be a normal subgroup of $G$. Let $\cala$ be an
additive category with a right $G$-action such that $N$ acts trivially.
Let $T$ be a left $G$-set such that $N$ operates freely.
Then there is an additive functor which is an equivalence of categories
\[
\cala \ast_G T \xrightarrow{\simeq} \cala \ast_{G/N} (N \backslash T).
\]
\item \label{ast-three}
Let $H$ be a subgroup of $G$ and $\cala$ be an additive category with
right $G$-action.
We denote by $\res_H \cala$ the additive category $\cala$ considered
with the $H$ action obtained by restriction.
Then for an $H$-set $T$ the map $T \to G \x_H T$ defined by
$t \mapsto (1_G,t)$ induces an equivalence of
additive categories
\[
 ( \res_H \cala ) \ast_H T \to \cala \ast_G (G \x_H T).
\]
\end{enumerate}
\end{proposition}

\begin{proof}
\ref{ast-first} The functor
$F \colon ( \cala \ast_K S ) \ast_G T \to
   \cala \ast_{K \times G} ( S \times T )$ is given
by
\[
\left( F(A) \right)_{(s,t)} = (A_t)_s \quad \mbox{ and } \quad
\left( F( \phi ) \right)_{(k,g)(s,t)} = ( \phi_{g,t} )_{k,s}.
\]
Note that if $\phi \colon A \to B$ is a morphism
in $( \cala \ast_K S ) \ast_G T$ then
$\phi_{g,t} \colon A_t \to g^{\ast}( B_{gt} )$ is a morphism in
$\cala \ast_K S$ and
$( \phi_{g,t} )_{k,s} \colon
 (A_t)_s \to k^{\ast} (( g^{\ast} B_{gt} )_{ks} )$
is a morphism in $\cala$. The target of this last morphism
is
\[
k^{\ast} ((g^{\ast}  B_{gt} )_{ks} ) =
   k^{\ast} ( (B_{gt} )_{ksg^{-1}} ) =
   k^{\ast} ( F(B)_{(ksg^{-1} ,gt )} ) =
(k,g)^{\ast} ( F(B)_{(k,g)(s,t)} )
\]
by definition of the $G$-action on $\cala \ast_K S$ and
the $K \times G$-action on $\cala$ and $S \times T$.
In particular this is a
correct target for $F( \phi )_{(k,g)(s,t)}$.
Clearly, $F$ is an isomorphism of categories.
To verify that $F$ is indeed  an  additive functor is
lengthy but straightforward.
\\
\ref{ast-second}
Let $p \colon T \to N \backslash T$ denote the projection
and consider
$\cala \ast_G p \colon \cala \ast_G T \to
 \cala \ast_G (N \backslash T)$,
see Remark~\ref{rem:cala-ast-T-is-natural}.
Next we define an additive functor
$F \colon \cala \ast_G (N \backslash T) \to
   \cala \ast_{G/N} (N \backslash T)$.
For both these categories objects are given by
sequences $(A_t)_{t \in T}$ of objects in $\cala$ indexed by $T$
and we define $F$ as the identity on objects.
Let $\phi \colon A \to B$ be a morphism in
$\cala \ast_G (N \backslash T)$.
For $g \in G$, $n \in N$ and $t \in T$
$\phi_{gn,Nt}$ is a morphism
\[
A_{Nt} \to (gn)^\ast B_{gnNt} = g^\ast B_{gNt}
\]
and we define $F$ on morphisms by
\[
(F(\phi))_{gN,Nt} = \sum_{n \in N} \phi_{gn,Nt}.
\]
Then the composition $F \circ (\cala \ast_G p)$
is full and faithful and hence an equivalence of
additive categories.
\\
\ref{ast-three}
It is straight forward to check that this functor is
full and faithful.
\end{proof}

Let $\Phi \colon K \to G$ be a group
homomorphism.
For a given additive category $\cala$ with right $K$-action
we define $\ind_{\Phi} \cala$,
a category with right $G$-action, as
\[
\ind_{\Phi} \cala = \cala \ast_K \res_{\Phi} G.
\]
This is a special case of the construction discussed before
Proposition~\ref{prop:main-technical-ast-properties}.
Here $\res_{\Phi} G$ denotes $G$
considered as $K$-$G$ biset via $\Phi$. Our main motivation for
Proposition~\ref{prop:main-technical-ast-properties} was the
following corollary which will play a key role when
we will study inheritance properties for Isomorphism Conjectures.

\begin{corollary} \label{cor:ind-res}
Let $\cala$ be an additive category with a right $K$-action.
For a group homomorphism $\Phi \colon K \to G$ and a $G$-set $T$
there is an additive functor which is an equivalence of categories
\[
( \ind_{\Phi} \cala ) \ast_{G} T \xrightarrow{\simeq} \cala \ast_{K}
      ( \res_{\Phi} T ).
\]
\end{corollary}

\begin{proof}
According to
Proposition~\ref{prop:main-technical-ast-properties}\ref{ast-first}
and
\ref{ast-second} we have an isomorphism respectively
an equivalence
\[
(\cala \ast_K \res_{\Phi} G ) \ast_G T \xrightarrow{\cong}
\cala \ast_{K \times G} (G \times T ) \xrightarrow{\simeq}
\cala \ast_K G \times_G T = \cala \ast_K \res_{\Phi} T.
\]
\end{proof}

\section{The Farrell-Jones  Conjecture with coefficients}
\label{sec:FJ-with-coeff}

Let $\IK^{-\infty} \colon \AddCat \to \Sp$ be the functor that
associates the non-connective $K$-theory spectrum to an
additive category (using the split exact structure).
This functor is constructed in
\cite{Pedersen-Weibel-delooping}.
See
\cite[Section 2.1 and 2.5]{BFJR-TOP}
for a brief review of this functor and its properties.
Let $G$ be a group and $\Or G$ be the orbit category of $G$,
whose objects are transitive $G$-sets of the form $G/H$
and whose morphisms are $G$-equivariant maps.
For any $\Or G$-spectrum $\bfE$, i.e.\
for any functor $\bfE \colon \Or G \to \Sp$, Davis and L\"uck
construct a $G$-equivariant homology theory for $G$-$CW$-complexes by
\[
H^G_*( X ; \bfE) = \pi_*( X_+ \sma_{\Or G} \bfE),
\]
where $X_+ \sma_{\Or G} \bfE$ denotes the balanced smash product of
$X_+ = \map_G ( ? , X_+)$ considered as a contravariant
$\Or G$-space and the covariant $\Or G$-spectrum $\bfE$.
For more details see
\cite[Section 4]{Davis-Lueck-assembly}.
For a group $G$ and a family $\calf$ of subgroups,
i.e.\ a collection of subgroups that is closed under subconjugation,
there is a $G$-$CW$-complex $E_\calf G$ with the property that for a subgroup
$H$ of $G$ the set of fixed points $E_\calf G^H$ is empty if $H \notin \calf$
and contractible if $H \in \calf$, see for example
\cite{Lueck-survey-classifying-spaces}.
The triple $(\bfE, \calf, G)$ is said to satisfy
the Isomorphism Conjecture if the so called assembly map
\[
H^G_*(E_\calf G;\bfE) \to H^G_*(\pt;\bfE) = \pi_*(\bfE(G/G))
\]
induced by the projection $E_\calf G \to \pt$ is an isomorphism,
see
\cite[Definition 5.1]{Davis-Lueck-assembly}.
In this paper we will use the following $\Or G$-spectra.

\begin{definition} \label{def:KuntenA}
Let $\cala$ be an additive category with right $G$-action.
The $\Or G$-spectrum $\bfK_\cala$ is defined by
\[
\bfK_\cala (T) = \IK^{-\infty}(\cala \ast_G T).
\]
\end{definition}

\begin{conjecture}[Algebraic $K$-Theory Farrell-Jones-Conjecture with Coefficients]
\label{con:FJC-with-coefficients-alg-K-theory}
Let $G$ be a group and let $\VCyc$ be the family of virtually cyclic subgroups of $G$.
Let $\cala$ be an additive category with a right $G$-action.
Then the assembly map
\[
H^G_*(E_\VCyc G;\bfK_\cala) \to H^G_*(\pt;\bfK_\cala)
\]
is an isomorphism.
\end{conjecture}

\begin{remark}
\label{rem:with-implies-without}
If $\cala = R_\oplus$ then by
Example~\ref{ex:cala-ast-T-and-Davis-Lueck}
$\bfK_\cala$ can be identified with the functor introduced in
\cite[Section 2]{Davis-Lueck-assembly}.
In particular Conjecture~\ref{con:FJC-with-coefficients-alg-K-theory}
implies the original conjecture by Farrell and Jones in
\cite{FJ-iso-conj}.
(The formulation of Davis and L\"uck has been
identified with the original formulation of
Farrell and Jones in
\cite{Hambleton-Pedersen-identify}.)
\end{remark}

\section{Inheritance properties}
\label{sec:inheritance-properties}

By definition a family of subgroups of a group $G$ is a collection of subgroups
closed under taking subgroups and conjugation.
If $\Phi \colon K \to G$ is group homomorphism and $\calf$ is a
family of subgroups of $G$ then we define a family
of subgroups of $K$ by setting
\[
\Phi^{\ast} \calf =
\{ H \subset G \; | \; H \mbox{ is a subgroup of } K
                \mbox{ and } \Phi ( H ) \in \calf \}.
\]

\begin{remark}
The $K$-$CW$-complex $\res_{\Phi} E_\calf G  $ is a model for
the classifying space $E_{  \Phi^{\ast} \calf } K$,
because it satisfies
the characterizing property concerning the
fixed point sets.
\end{remark}

\begin{proposition}
\label{prop:equivalence-of-assembly-maps}
Let $\Phi \colon K \to G$ be a group homomorphism.
Let $\cala$ be an additive category with right $K$-action and let
$\calf$ be a family of subgoups of $G$. Then
the assembly map
\[
H^G_*( E_\calf G ; \bfK_{\ind_{\Phi} \cala}) \to
H^G_*( \pt; \bfK_{\ind_{\Phi} \cala})
\]
is equivalent to the assembly map
\[
H^K_*( E_{ \Phi^{\ast}\calf } K; \bfK_\cala ) \to
H^K_*( \pt; \bfK_\cala ).
\]
\end{proposition}

\begin{proof}
Because of Corollary~\ref{cor:ind-res} and
since $\IK^{-\infty}$
preserves equivalences we have equivalences
of $\Or G$-spectra
$\bfK_{ \ind_{\Phi} \cala} \xrightarrow{\simeq}
 \bfK_\cala \circ \res_\Phi$
and therefore for every $G$-space $X$ a natural
isomorphism
\[
H^G_*(X; \bfK_{\ind_\Phi \cala}) \cong
H^G_*(X; \bfK_{\cala} \circ \res_\phi)
\]
For a $K$-space $Y$ define $\ind_{\Phi} Y$ to be the quotient
of $G \times Y$ by the right $K$ action given by
$( g,y)k = (g \Phi(k) , k^{-1} y)$.
For every $G$-space $X$ there is an isomorphism
\begin{eqnarray*}
\res_{\Phi} X^? _+ & = &  \map_K ( ? , \res_{\Phi} X )_+ \\
& \cong & \map_G ( \ind_{\Phi} (?) , X )_+ \\
& \cong & \map_G ( ?? , X )_+ \sma_{\Or G} \map_G ( \ind_{\Phi} (?) , ?? )_+ \\
& = & X_+^{??} \sma_{\Or G} \map_G ( \ind_{\Phi} (?) , ?? )_+
\end{eqnarray*}
of contravariant pointed $\Or K$-spaces.
Here and in the next display $?$ denotes functoriality in $\Or K$
and $??$ denotes functoriality in $\Or G$.
For every covariant functor $\bfF$ from $K$-sets to spectra there is
an isomorphism of covariant $\Or G$-spectra
\begin{eqnarray*}
\map_G ( \ind_{\Phi} ( ? ) , ?? )_+ \sma_{\Or K } \bfF ( ? ) &
\cong & \map_K ( ? , \res_{\Phi} ( ?? ))_+ \sma_{\Or K} \bfF ( ? ) \\
& \cong  & \bfF \circ \res_{\Phi} ( ?? ) .
\end{eqnarray*}
Combining these isomorphisms with associativity of
balanced smash products we obtain an isomorphism of spectra
\[
\res_{\Phi} X_+ \sma_{\Or K } \bfF \cong X_+
\sma_{ \Or G } \bfF \circ \res_{\Phi}
\]
and in particular a natural isomorphism
\[
H^K(\res_\phi X; \bfK_\cala) \cong H^G(X; \bfK_\cala \circ \res_\Phi).
\]
It remains to observe that $\res_{\Phi} \pt = \pt $ and that
$\res_{\Phi} E_{ \calf } G$  is a model for $E_{  \Phi^{\ast} \calf } K$.
\end{proof}

Proposition~\ref{prop:equivalence-of-assembly-maps} has
the following immediate consequence.

\begin{corollary}
\label{cor:Inheritance-for-group-homoeomorphisms}
Let $\Phi \colon K \to G$ be a group homomorphism.
Let $\calf$ be a family of subgroups of $G$.
Suppose that for every additive category $\cala$ with $G$-action
the assembly map
\[
H_*^G(E_\calf G;\bfK_\cala) \to H_*^G(\pt;\bfK_\cala)
\]
is injective.
Then for every additive category $\calb$ with $K$-action
the assembly map
\[
H_*^K(E_{\Phi^*\calf} K;\bfK_\calb) \to H_*^K(\pt;\bfK_\calb)
\]
is injective. The same statement holds with injectivity
replaced by surjectivity in assumption and conclusion.
\end{corollary}

\begin{remark}[With coefficients is stronger than fibered]
\label{rem:coefficients-implies-fibered}
The fibered version of the Farrell-Jones Conjecture
in algebraic $K$-theory for a group $G$ (and a ring $R$),
\cite[Section 1.7]{FJ-iso-conj}
can be formulated as follows:
for every group homeomorphism $\Phi \colon K \to G$
the assembly map
$H^K_*(E_{\Phi^*\VCyc} K;\bfK_R) \to H^K_*(\pt;\bfK_R)$
is an isomorphism, see Section~6 and in particular Remark~6.6 in
\cite{Bartels-Lueck-trees}.
Therefore by Corollary~\ref{cor:Inheritance-for-group-homoeomorphisms}
the Farrell-Jones Conjecture with
Coefficients~\ref{con:FJC-with-coefficients-alg-K-theory}
implies the Fibered Farrell-Jones Conjecture.
\end{remark}

Corollary~\ref{cor:Inheritance-for-group-homoeomorphisms}
implies in particular the following theorem about
the Farrell-Jones Conjecture with
Coefficients~\ref{con:FJC-with-coefficients-alg-K-theory}.

\begin{theorem} \label{thm-passage-to-subgroups}
\label{thm:FJ-w-C-passes-to-subgroups-alg-K}
Let $H$ be a subgroup of $G$.
Suppose that for every additive category $\cala$ with $G$-action
the assembly map
\[
H_*^G(E_\VCyc G;\bfK_\cala) \to H_*^G(\pt;\bfK_\cala)
\]
is injective or surjective respectively.
Then for every additive category $\calb$ with right $H$-action
the assembly map
\[
H_*^H(E_\VCyc H;\bfK_\calb) \to H_*^H(\pt;\bfK_\calb)
\]
is injective or surjective respectively.
In particular, if the
Farrell-Jones Conjecture with
Coefficients~\ref{con:FJC-with-coefficients-alg-K-theory}
holds for a group $G$,
then it holds for every subgroup of $G$.
\end{theorem}

Similar as for rings there exists a suspension category and hence results that hold without a condition
on the coefficient category can always be shifted down. More precisely the following holds.

\begin{proposition}
\label{prop:delooping-of-assembly-map}
For every additive category $\cala$ with $G$-action
there is an additive category $\Sigma \cala$ with
$G$-action such that for
every family of subgroups $\calf$ and
every $n \in \IZ$ the assembly map
\begin{eqnarray*}
H^G_n ( E_\calf G;\bfK_\cala)
 \to H^G_n(\pt;\bfK_\cala)
\end{eqnarray*}
is isomorphic to the assembly map
\begin{eqnarray*}
H^G_{n-1} ( E_\calf G;\bfK_{\Sigma \cala})
 \to H^G_{n-1}(\pt;\bfK_{\Sigma \cala}).
\end{eqnarray*}

\end{proposition}

\begin{proof}
We use a construction of $\Sigma \cala$
that is similar to the construction
from \cite{Pedersen-Weibel-delooping}.
For a given $\cala$ there is a natural construction of
a Karoubi filtration of additive categories
$\cala' \subset \Lambda \cala$ whose quotient we denote by
$\Sigma \cala$. Here $\cala'$ is naturally
equivalent to $\cala$ and there is an Eilenberg swindle on
$\Lambda \cala$, see Example~\ref{ex:lambda-cala}.
Therefore $\cala \to \Lambda \cala \to \Sigma \cala$
induces a fibration sequence in (non-connective) $K$-theory by
\cite[Theorem 1.28]{Carlsson-Pedersen-Controlled-algebra-Novikov}
and $K_* \Lambda \cala = 0$.
Because the construction is natural, there are $G$-actions
on $\Lambda \cala$ and $\Sigma \cala$.
Both,  the Karoubi filtration and the Eilenberg swindle are
preserved by the passage from $\cala$ to  $\cala \ast_G T$.
Therefore we have a fibration sequence of $\Or G$-spectra,
\[
\bfK_\cala \to \bfK_{\Lambda \cala} \to \bfK_{\Sigma \cala}
\]
that gives long exact sequences of the associated
homology groups for every $G$-space $X$.
By the Eilenberg swindle on $\cala \ast_G T$ the groups
$H^G_*(X;\bfK_{\Lambda \cala})$ vanish and the boundary map
in the long exact sequence yields the desired
identification of assembly maps.
\end{proof}

From Proposition~\ref{prop:delooping-of-assembly-map}
we obtain the following analog of
\cite[Corollary 7.3]{BFJR-TOP}.

\begin{corollary}
\label{cor:inheritance-n-to-n-1}
Let $\calf$ be a family of subgroups of the group $G$.
If for every additive category $\cala$ with right $G$-action
the assembly map
\[
H^G_* ( E_\calf G;\bfK_\cala)  \to H^G_*(\pt;\bfK_\cala)
\]
is injective or surjective respectively in a fixed degree
$* = n$, then this assembly map is
injective or surjective respectively in all degrees $* = j$
with $j \leq n$.
\end{corollary}

\section{$L$-theory}
\label{sec:L-theory}

Everything we did for algebraic $K$-theory
has an analog in $L$-theory and we will state
the corresponding conjecture and inheritance result
here quickly.
An additive category with involution is an additive category $\cala$
together with an additive contravariant functor
$\# = ( - )^{\#} \colon \cala \to \cala$
such that $\# \circ \# = \id$.
We consider now additive categories
with involution and right $G$-action, where we
require in addition that for every $g \in G$
the covariant functor $g^{\ast}$
is compatible with the involution $\#$,
i.e.\ $\# \circ g^{\ast} = g^{\ast} \circ \#$.
If $T$ is a $G$-set then
\[
(A^\#)_t = (A_t)^\#  \quad \mbox{and} \quad
(\phi^\#)_{g,t} = g^\ast ( (\phi_{g^{-1},gt} )^\#)
\]
defines an involution on $\cala *_G T$.
There is a functor $\IL^{-\infty} \colon \AddCatInv \to \Sp$ that
associates the $L$-theory spectrum
to an additive category with involution constructed by
Ranicki
\cite{Ranicki-L-theory-top-manifolds}.
We consider the $\Or G$-spectrum
$\bfL_\cala$ defined by
\[
\bfL_\cala (T) = \IL^{-\infty}(\cala \ast_G T).
\]

\begin{conjecture}[$L$-Theory Farrell-Jones-Conjecture
                   with Coefficients]
\label{con:FJC-with-coefficients-L-theory}
Let $G$ be a group and let $\VCyc$ be the family of virtually cyclic
subgroups of $G$.
Let $\cala$ be an additive category with involution
with a right $G$-action.
Then the assembly map
\[
H^G(E_\VCyc G);\bfL_\cala) \to H^G(\pt;\bfL_\cala)
\]
is an isomorphism.
\end{conjecture}

The only property of the functor $\IK^{-\infty}$
that was used in the proof of
Proposition~\ref{prop:equivalence-of-assembly-maps}
is that it sends equivalences of categories to
equivalences of spectra.
Because this
property holds also for the functor $\IL^{-\infty}$
there is also the $L$-theory version of
Proposition~\ref{prop:equivalence-of-assembly-maps}.
Therefore there are also $L$-theory versions of
Corollary~\ref{cor:Inheritance-for-group-homoeomorphisms},
Remark~\ref{rem:coefficients-implies-fibered} and
Theorem~\ref{thm:FJ-w-C-passes-to-subgroups-alg-K}.
We spell out only the analog of Theorem~\ref{thm:FJ-w-C-passes-to-subgroups-alg-K}.

\begin{theorem}
\label{thm:FJ-w-C-passes-to-subgroups-L}
Let $H$ be a subgroup of $G$.
Suppose that for every additive category $\cala$
with involution with $G$-action
the assembly map
\[
H_*^G(E_\VCyc G;\bfL_\cala) \to H_*^G(\pt;\bfL_\cala)
\]
is injective or surjective respectively.
Then for every additive category $\calb$
with involution with right $H$-action
the assembly map
\[
H_*^H(E_\VCyc H;\bfL_\calb) \to H_*^H(\pt;\bfL_\calb)
\]
is injective or surjective respectively.
In particular, if the $L$-Theory
Farrell-Jones Conjecture with
Coefficients~\ref{con:FJC-with-coefficients-L-theory}
holds for a group $G$,
then it also holds for every subgroup of $G$.
\end{theorem}

\section{Crossed products} \label{sec-crossed}
\label{sec:crossed-products}

In this section we show that the
Farrell-Jones Conjecture with
Coefficients~\ref{con:FJC-with-coefficients-alg-K-theory}
covers crossed product rings. We first recall the notion of a crossed product
ring, compare \cite{Passman(1989)}.

Let $R$ be a ring, $G$ be a group and
$\alpha \colon G \to \Aut(R), g \mapsto \alpha_g$ and
$\tau \colon G \x G \to R^{\x}, (g,h) \mapsto \tau_{g,h}$ be maps.
Here $R^\x$ are the units of $R$ and $\Aut(R)$ denotes the group of ring-automorphisms of $R$.
We require that
\begin{eqnarray}
\label{eq:crossed-tau-tau}
\tau_{g,h} \tau_{gh,k} & = & \alpha_g( \tau_{h,k}) \tau_{g,hk} \\
\label{eq:crossed-alpha-alpha}
\tau_{g,h} \alpha_{gh}(r) & = & (\alpha_g \circ \alpha_h)(r) \tau_{g,h}
\end{eqnarray}
for $g,h,k \in G$, $r \in R$.
We will also assume that $\alpha_e = \id_R$, where $e$ denotes
the unit element in $G$.

The crossed product ring $R_{\alpha,\tau} G$ is as
an additive group $RG$, but is equipped with a twisted multiplication
$\cdot_{\alpha,\tau}$ where
\begin{equation}
\label{eq:crossed-twisted-product}
(r g) \cdot_{\alpha,\tau} (s h) = r \alpha_g(s) \tau_{g,h} g h
\end{equation}
for $r,s \in R$ and $g,h \in G$.
The element $1_R e$ is the unit of $R_{\alpha,\tau} G$.

Setting $g = e$ or $h = e$ in (\ref{eq:crossed-alpha-alpha})
we conclude that for $g \in G$,
\begin{equation}
\label{eq:tau-center}
\tau_{e,g} \; \mbox{and} \; \tau_{g,e} \;
  \mbox{lie in the center of} \; R.
\end{equation}

\begin{example}
The  notion of a crossed product ring naturally appears in the following situation.
Consider an extension of groups
\[
1 \to K \to \Gamma \xrightarrow{p} G \to 1.
\]
Choose a set-theoretical section $s$ of $p$ such that $s(1) = 1$. Let $S$ be a ring, let $R=SK$ and
set $\alpha_g (r) = s(g) r s(g)^{-1}$ and  $\tau_{g,h} = s(g) s(h ) s ( gh )^{-1}$.
Then
\[
R_{\alpha , \tau} G \cong S \Gamma.
\]
\end{example}
Our aim is now to define an additive category
$\cala_{\alpha,\tau}$
with a right $G$-action such that the category $\cala \ast_G \pt$ is equivalent
to the category of finitely generated free
$R_{\alpha,\tau} G$-modules.

We start with the category $\calf^f(R)$ of finitely
generated free $R$-modules.
If $\varphi$ is an automorphism of $R$, then we define
a functor $M \mapsto \res_\varphi M$ where the latter is the
$R$-module obtained by twisting the $R$-module structure by
$\varphi$, i.e.\
$r \cdot_{\res_\varphi M} v = \varphi(r) \cdot_M v$.
This defines a right action of $\Aut(R)$ on $\calf^f(R)$.
(If we want a small category we can restrict attention
to modules of the form $\res_\varphi R^n$.)

We digress for a moment and discuss the special case where
$\tau \equiv 1_R$.
Then $\alpha$ is a group homomorphism
and we obtain an action of $G$ on $\calf^f(R)$.
This is the desired category
with $G$-action in this special case.
The equivalence to the category of finitely generated
free $R_{\alpha} G$-modules sends
the morphism $\phi \colon M \to N$ in
$\calf^f (R) \ast_G \pt$ with components
\[
\phi_g \colon M \to \res_{\alpha_g} N
\]
to the $R_{\alpha}G$-linear map
\[
R_{\alpha} G \otimes_R  M \to R_{\alpha} G \otimes_R N,
         \quad x \otimes v \mapsto xg^{-1} \otimes \phi_g (v).
\]

We continue with the explanation of  the general case.
In general, $L_{\tau_{g,h}}$ defines a
natural transformation from
$\res_{\alpha_{gh}}$ to $\res_{\alpha_h} \circ \res_{\alpha_g}$.
(Here we denoted the map $v \mapsto rv$ for $r \in R$ by $L_r$.
The expression $rv$ is formed with respect to the the
original module multiplication on $M$.)

The category $\cala_{\alpha,\tau}$ is now obtained by rigidifying
this operation as follows.
Objects of $\cala_{\alpha,\tau}$ are pairs $(M,g)$ where
$M$ is a finitely generated free $R$-module and
$g \in G$.
Morphisms from $(M,g)$ to $(N,h)$ are $R$-linear maps
$\varphi \colon \res_{\alpha_g} M \to \res_{\alpha_h} N$
and composition is composition of linear maps.
The right action of $\gamma \in G$ is defined by
$(M,g) \mapsto (M,g\gamma)$ on objects and by
\begin{equation}
\label{eq:gamma-varphi}
\varphi \mapsto \gamma^*\varphi = L_{\tau_{h,\gamma}}^{-1} \circ
      \varphi \circ L_{\tau_{g,\gamma}}
\end{equation}
for a morphism $\varphi \colon (M,g) \to (N,h)$.
The only thing one has to check is that
$\delta^* (\gamma^* \varphi) = (\gamma \delta)^*\varphi$
for a morphism $\varphi \colon (M,g) \to (N,h)$.
Recall that such a morphism is given by
an additive map $\varphi \colon M \to N$ for which
$\varphi \circ L_{\alpha_g(r)} = L_{\alpha_h (r)} \circ \varphi$
for $v \in M$, $r \in R$.
Using this and \eqref{eq:crossed-tau-tau}
we can compute
\begin{eqnarray*}
\delta^* (\gamma^*\varphi) & = &
  L_{\tau_{h\gamma,\delta}}^{-1} \circ L_{\tau_{h,\gamma}}^{-1}
    \circ \varphi \circ L_{\tau_{g,\gamma}}
    \circ L_{\tau_{g \gamma,\delta}} \\
  & = &
  L_{\tau_{h, \gamma\delta}}^{-1} \circ L_{\alpha_h(\tau_{\gamma,\delta})}^{-1}
    \circ \varphi \circ L_{\alpha_g(\tau_{\gamma,\delta})}
    \circ L_{\tau_{g, \gamma \delta}} \\
  & = &
  L_{\tau_{h, \gamma\delta}}^{-1} \circ L_{\alpha_h(\tau_{\gamma,\delta})}^{-1}
    \circ L_{\alpha_h(\tau_{\gamma,\delta})} \circ \varphi
    \circ L_{\tau_{g, \gamma \delta}} \\
  & = &
  L_{\tau_{h, \gamma\delta}}^{-1}
    \circ \varphi
    \circ L_{\tau_{g, \gamma \delta}} \\
  & = &
  (\gamma \delta)^* \varphi.
\end{eqnarray*}

\begin{proposition}
\label{prop:crossed-products-and-categories}
The categories $\cala_{\alpha,\tau} *_G \pt$ and the
category $\calf^f (R_{\alpha,\tau} G)$ of finitely generated free
$R_{\alpha,\tau} G$ modules are equivalent as additive categories.
\end{proposition}

\begin{proof}
We start by listing a number of useful consequences of
(\ref{eq:crossed-tau-tau}), (\ref{eq:crossed-alpha-alpha}) and
$\alpha_e = \id_R$,
\begin{eqnarray}
\alpha_a(\tau_{b,c}) & = & \tau_{a,b} \tau_{ab,c} \tau_{a,bc}^{-1}
      \label{eq:alpha-tau}
\\
\alpha^{-1}_a (r) & = & \tau^{-1}_{a^{-1},a} \alpha_{a^{-1}}(r)
      \tau_{a^{-1},a}
      \label{eq:alpha-inv-tau}
\\
(\alpha_{a} \circ \alpha_b) (r) & = &
      \tau_{a,b} \alpha_{ab}(r) \tau_{a,b}^{-1}
      \label{eq:alpha-alpha-r}
\end{eqnarray}
for $a$, $b$, $c \in G$ and $r \in R$.
From the definition of the product $\cdot = \cdot_{\alpha,\tau}$
in (\ref{eq:crossed-twisted-product}) we recall
\begin{eqnarray}
a \cdot r & = & \alpha_a(r) \cdot a
         \label{eq:a-cdot-r}
\\
r \cdot a & = & a \cdot \alpha_a^{-1}(r)
         \label{eq:r-cdot-a}
\\
a \cdot b & = & \tau_{a,b} \cdot ab
         \label{eq:a-cdot-b}
\end{eqnarray}
for $a,b \in G$ and $r \in R$.

Denote by $(\cala_{\alpha,\tau}*_G \pt )_0$ the full subcategory
of $\cala_{\alpha,\tau} *_G \pt$ whose objects are of the form
$(M,e)$.
It is easy to check that the inclusion
$(\cala_{\alpha,\tau} *_G \pt )_0 \to \cala_{\alpha,\tau} *_G \pt$
is an equivalence of additive categories.
We define a functor
$F \colon (\cala_{\alpha,\tau} *_G \pt)_0 \to \calf^f (R_{\alpha,\tau} G)$
as follows.
For an object $(M,e)$ in $(\cala_{\alpha,\tau} *_G \pt )_0$
let $F(M,e) = R_{\alpha,\tau} G \ox_R M$.
A morphism $\phi \colon (M,e) \to (N,e)$ in
$(\cala_{\alpha,\tau} *_G \pt)$
is by definition a sequence $(\phi_\gamma)_{ \gamma \in G}$
where $\phi_\gamma \colon M \to \res_{\alpha_\gamma} N$
is an $R$-linear map.
Because we can add morphisms in additive categories
it will suffice to discuss morphisms for which
$\phi_\gamma = 0$ for all but one $\gamma \in G$;
we write $(\varphi,g)$ for the morphism given by
$\phi_\gamma = \varphi$ if $\gamma = g$ and $\phi_\gamma = 0$
otherwise, in particular $\varphi$ is an additive map $M \to N$
for which
\begin{equation}
\label{eq:varphi-r}
\varphi ( r v ) = \alpha_g(r) \varphi (v) \quad \mbox{for all}
     \quad r \in R, v \in M.
\end{equation}
We define $F(\varphi,g) \colon F(M,e) \to F(N,e)$
as the linear map
\begin{equation}
\label{eq:define-F(varphi,g)}
x \ox v \mapsto x \cdot g^{-1} \cdot \tau^{-1}_{g,g^{-1}} \ox \varphi(v).
\end{equation}
Note that
\begin{align*}
 xr \cdot g^{-1} \cdot \tau^{-1}_{g,g^{-1}}
 & = x \cdot g^{-1} \cdot \alpha^{-1}_{g^{-1}}(r) \tau^{-1}_{g,g^{-1}}
 & \mbox{by} \; (\ref{eq:r-cdot-a})
\\
 & = x \cdot g^{-1} \cdot \tau^{-1}_{g,g^{-1}} \alpha_g(r)
 & \mbox{by} \; (\ref{eq:alpha-inv-tau}) &.
\end{align*}
Because of (\ref{eq:varphi-r}) this means
that $F(\varphi,g)$ defines indeed a well defined map
on the tensor product.
(This explains the appearance of $\tau_{g,g^{-1}}^{-1}$ in
(\ref{eq:define-F(varphi,g)}); without this term the map
is ill defined.)

Next we check that
$F$ is compatible with composition, a somewhat
tedious calculation.
Let $(\psi,h) \colon (N,e) \to (L,e)$ be a second morphism
in $(\cala_{\alpha,\tau} *_G \pt )_0$,
in particular $\psi$ is an additive map $N \to L$
for which
\begin{equation}
\label{eq:psi-r}
\psi ( r v ) = \alpha_h(r) \psi (v) \quad \mbox{for all}
     \quad r \in R, v \in N.
\end{equation}
Then $F(\psi,h) \circ F(\varphi,g)$ maps $x \ox v$ to
\begin{align*}
&     x \cdot g^{-1} \cdot \tau_{g,g^{-1}}^{-1} \cdot h^{-1} \cdot
          \tau^{-1}_{h,h^{-1}} \ox \psi(\varphi(v))
\\
& =   x \cdot g^{-1} \cdot h^{-1} \cdot
       \alpha^{-1}_{h^{-1}}(\tau^{-1}_{g,g^{-1}})
       \tau^{-1}_{h,h^{-1}} \ox \psi(\varphi(v))
&                \mbox{by} \; (\ref{eq:r-cdot-a})
\\
& =   x \cdot \tau_{g^{-1},h^{-1}} \cdot (hg)^{-1} \cdot
       \alpha^{-1}_{h^{-1}}(\tau^{-1}_{g,g^{-1}})
       \tau^{-1}_{h,h^{-1}} \ox \psi(\varphi(v))
&                \mbox{by} \; (\ref{eq:a-cdot-b})
\\
& =   x \cdot (hg)^{-1} \cdot \alpha^{-1}_{(hg)^{-1}}(\tau_{g^{-1},h^{-1}})
       \alpha^{-1}_{h^{-1}} (\tau^{-1}_{g,g^{-1}})
       \tau^{-1}_{h,h^{-1}} \ox \psi(\varphi(v))
&                \mbox{by} \; (\ref{eq:r-cdot-a})
\\
& =   x \cdot (hg)^{-1} \cdot A \ox \psi(\varphi(v))
\end{align*}
where
\begin{align*}
A  = \; & \alpha^{-1}_{(hg)^{-1}} (\tau_{g^{-1},h^{-1}})
        \alpha^{-1}_{h^{-1}} (\tau^{-1}_{g,g^{-1}})
        \tau^{-1}_{h,h^{-1}}
\\
   = \; & \left(
            \tau^{-1}_{hg,(hg)^{-1}}
            \alpha_{hg} (\tau_{g^{-1},h^{-1}})
            \tau_{hg,(hg)^{-1}}
        \right)
\\
        & \left(
            \tau^{-1}_{h,h^{-1}}
            \alpha_h (\tau^{-1}_{g,g^{-1}})
            \tau_{h,h^{-1}}
        \right)
        \tau_{h,h^{-1}}^{-1}
      & \mbox{by} \; (\ref{eq:alpha-inv-tau})
\\
   = \; & \tau^{-1}_{hg,(hg)^{-1}}
        \left(
        \tau^{-1}_{h,g}
        \alpha_h ( \alpha_g (\tau_{g^{-1},h^{-1}} ) )
        \tau_{h,g}
        \right)
\\
      & \tau_{hg,(hg)^{-1}}
        \tau^{-1}_{h,h^{-1}}
        \alpha_h (\tau^{-1}_{g,g^{-1}})
      & \mbox{by} \; (\ref{eq:alpha-alpha-r})
\\
   = \; & \tau^{-1}_{hg,(hg)^{-1}}
        \tau^{-1}_{h,g}
        \left(
        \alpha_h ( \tau_{g,g^{-1}}
                   \tau_{e,h^{-1}}
                   \tau_{g,(hg)^{-1}}^{-1}
                 )
        \right)
\\
      & \tau_{h,g}
        \tau_{hg,(hg)^{-1}}
        \tau^{-1}_{h,h^{-1}}
        \alpha_h (\tau^{-1}_{g,g^{-1}})
      & \mbox{by} \; (\ref{eq:alpha-tau})
\\
   = \; & \tau^{-1}_{hg,(hg)^{-1}}
        \tau^{-1}_{h,g}
        \left(
        \alpha_h ( \tau_{g,g^{-1}} )
        \alpha_h ( \tau_{e,h^{-1}} )
        \alpha_h ( \tau_{g,(hg)^{-1}} )^{-1}
        \right)
        \\
        &
        \tau_{h,g}
        \tau_{hg,(hg)^{-1}}
        \tau^{-1}_{h,h^{-1}}
        \alpha_h (\tau_{g,g^{-1}})^{-1}
\\
   = \; & \tau^{-1}_{hg,(hg)^{-1}}
        \tau^{-1}_{h,g}
        \left( \tau_{h,g} \tau_{hg,g^{-1}} \tau_{h,e}^{-1} \right)
        \left( \tau_{h,e} \tau_{h,h^{-1}} \tau_{h,h^{-1}}^{-1} \right)
\\
      & \left( \tau_{h,g} \tau_{hg,(hg)^{-1}} \tau^{-1}_{h,h^{-1}}
                                                        \right)^{-1}
        \tau_{h,g}
        \tau_{hg,(hg)^{-1}}
        \tau^{-1}_{h,h^{-1}}
        \left( \tau_{h,g} \tau_{hg,g^{-1}} \tau_{h,e}^{-1} \right)^{-1}
       & \mbox{by} \; (\ref{eq:alpha-tau})
\\
   = \; & \tau^{-1}_{hg,(hg)^{-1}}
        \tau_{hg,g^{-1}}
        \tau_{h,e}
        \tau_{hg,g^{-1}}^{-1}
        \tau_{h,g}^{-1}
\\
   = \; & \tau^{-1}_{hg,(hg)^{-1}}
        \tau_{h,g}^{-1}
        \tau_{h,e}
      & \mbox{by} \; (\ref{eq:tau-center}) & .
\end{align*}
Compute the composition in $(\cala_{\alpha,\tau} *_G \pt)_0$
as follows
\begin{align*}
(\psi,h) \circ (\varphi,g) & = (g^*\psi \circ \varphi, hg) \quad
      & \mbox{by} \;  (\ref{eq:composition-in-cala-ast})
\\
& = (L^{-1}_{\tau_{h,g}} \circ \psi \circ
         L_{\tau_{e,g}} \circ \varphi, hg) \quad
      & \mbox{by} \; (\ref{eq:gamma-varphi}) &.
\end{align*}
Therefore $F((\psi,h)\circ(\varphi,g))$
maps $x \ox v$ to
\begin{align*}
& (x \cdot (hg)^{-1}  \cdot  \tau^{-1}_{hg,(hg)^{-1}} \ox
         \tau^{-1}_{h,g} \psi(\tau_{e,g} \varphi(v))
\\
 & =   x \cdot (hg)^{-1} \cdot \tau^{-1}_{hg,(hg)^{-1}} \ox
         \tau^{-1}_{h,g} \alpha_h(\tau_{e,g})
         \psi(\varphi(v)) \qquad
       \mbox{by} \; (\ref{eq:psi-r})
\\
& =   x \cdot (hg)^{-1} \ox B \psi(\varphi(v))
\end{align*}
where
\begin{align*}
B & = \tau^{-1}_{hg,(hg)^{-1}} \tau^{-1}_{h,g}
        \alpha_h(\tau_{e,g})
\\
  & = \tau^{-1}_{hg,(hg)^{-1}} \tau^{-1}_{h,g}
        \tau_{h,e} \tau_{h,g} \tau_{h,g}^{-1}
  & \mbox{by} \; (\ref{eq:alpha-tau})
\\
  & = \tau^{-1}_{hg,(hg)^{-1}} \tau^{-1}_{h,g}
        \tau_{h,e}.
\end{align*}
Thus $A = B$ and this shows
$F(\psi,h) \circ F(\varphi,g) = F((\psi,h) \circ (\varphi,g))$.
Thus $F$ is indeed a functor.
It is straight forward to check that $F$ is full and faithful, i.e.\
an equivalence of categories.
\end{proof}

The following sharpening of Proposition~\ref{prop:crossed-products-and-categories}
is obtained by formal arguments.

\begin{corollary}
\label{cor:crossed-product-in-category-language}
Suppose we are given a crossed product situation
\[
R, \quad  \alpha \colon G \to \Aut (R ) , \quad \tau \colon G \times G \to R^{\times}.
\]
Then there exists
an additive category $\cala_{\alpha , \tau}$ with a right $G$-action, such that for every orbit $G/H$
the category
\[
\cala_{\alpha, \tau} \ast_G G/H \quad \mbox{ and the category } \quad
\calf^f (R_{\alpha|, \tau|} H )
\]
of finitely generated $R_{\alpha|, \tau|} H$-modules
are equivalent.
Here $\alpha|$ and $\tau|$ denote the restriction of $\alpha$ and $\tau$ to $H$ respectively $H \times H$.
In particular there is for every $G/H$ and every $n \in \IZ$ an isomorphism
\[
K_n ( \cala_{\alpha , \tau} \ast_G G/H ) \cong K_n ( R_{\alpha| , \tau| } H ).
\]
\end{corollary}
\begin{proof}
We have a chain of equivalences
\begin{eqnarray*}
\cala_{\alpha, \tau} \ast_G G/H & \simeq & ( \res_H \cala_{\alpha , \tau} ) \ast_H \pt \\
 & \simeq & \cala_{\alpha| , \tau|} \ast_H \pt \\
 & \simeq & R_{\alpha| ,\tau|} H_{\oplus} .
\end{eqnarray*}
Here the first equivalence is a special case of
Proposition~\ref{prop:main-technical-ast-properties}~\ref{ast-three} and the last follows immediately
from the previous Proposition~\ref{prop:crossed-products-and-categories}. The second equivalence
is induced from the $H$-equivariant inclusion
$\cala_{\alpha| , \tau|} \to \res_H \cala_{\alpha , \tau}$ which sends $(M,h)$ to the same element
considered as an object of $\cala_{\alpha , \tau}$. This inclusion is clearly full and faithful and every
object $(M,g)$ in the target is isomorphic to $(\res_{\alpha_g} M , e)$. One then checks that in general an $H$-equivariant
equivalence $\cala \to \calb$ induces an equivalence $\cala \ast_H T \to \calb_H \ast T$ for every $G$-set $T$.
\end{proof}

Observe that in particular the $G$-equivariant homology
theory $H_{\ast}^G( - ; \bfK_{\cala_{\alpha , \tau}} )$ evaluated on an orbit $G/H$ is isomorphic to
$K_{\ast} ( R_{\alpha| , \tau| } H )$. The following special case of Conjecture~\ref{con:FJC-with-coefficients-alg-K-theory}
hence makes precise the idea that $K_{\ast}( R_{\alpha , \tau} G)$ should be assembled from the
pieces $K_{\ast} ( R_{\alpha| , \tau| } H )$, where $H$ ranges over the virtually cyclic subgroups of $G$.
\begin{conjecture}
\label{con:FJ-for-crossed-product}
Suppose $R_{\alpha , \tau} G$ is a crossed product ring, then the assembly map
\[
H^G_*(E_\VCyc G;\bfK_{\cala_{\alpha , \tau}} ) \to H^G_*(\pt;\bfK_{\cala_{\alpha , \tau}} )\cong K_{\ast} ( R_{\alpha , \tau } G )
\]
induced from $E_\VCyc G \to \pt$ is an isomorphism.
\end{conjecture}

\section{Controlled algebra}
\label{sec:controlled-algebra}

Many results on the Farrell-Jones conjecture (without coefficients)
use the concept of controlled algebra.
In this section we briefly indicate how the fundamental concepts of controlled algebra
extend from rings to additive categories with group actions.

The following generalizes the definitions in
\cite[Section~2]{BFJR-TOP}.

\begin{definition}
Let $\cala$ be an additive category with a right $G$-action and
let $X$ be a free $G$-space.
Define the additive category with right $G$-action
\[
\calc ( X ; \cala )
\]
as follows.
Objects are families $A = (A_x )_{x \in X}$ of objects in $\cala$ such that
$\supp A = \{ x \in X \; | \; A_x \neq 0 \}$ is locally finite.
A morphism $\phi \colon A \to B$ is a family
$( \phi_{y,x} )_{(y,x) \in X \times X}$,
where $\phi_{y,x} \colon A_x \to B_y$
is a morphism in $\cala$ and for fixed $x$ the set of $y$ with
$\phi_{y,x} \neq 0$ is finite and for fixed $y$ the set of $x$ with
$\phi_{y,x} \neq 0$ is finite.
The composition $\psi = \phi^{\prime} \circ \phi$ is defined to be
\[
\psi_{z,x} = \sum_{y \in X} \phi^{\prime}_{z,y} \circ \phi_{y,x}.
\]
The element $g \in G$ acts via the covariant additive
functor $g^{\ast}$ which is given by
\[
(g^{\ast} A)_x = g^{\ast}(A_{gx} ) \quad \mbox{ and }
(g^{\ast} \phi)_{y,x} =  g^{\ast} ( \phi_{gy,gx} ).
\]
\end{definition}

It now makes sense to consider the fixed category $\calc( X ; \cala )^G$.
An object $A$  and a morphism $\phi$  in the fixed category satisfy
\[
A_x = g^{\ast} ( A_{gx} ) \quad \mbox{ and } \quad \phi_{y,x} =
g^{\ast} ( \phi_{gy, gx} ).
\]
Observe that in the case where the category $\cala$ is $R_{\oplus}$
for some ring $R$ equipped with the trivial $G$-action
we obtain the category which was denoted $\calc^{G} ( X ; R )$ in
\cite[Section~2]{BFJR-TOP}, compare also Example~\ref{ex:cala-ast-T-and-Davis-Lueck}.

\subsection{Support conditions}
As usual we define the support of an object $A$,
respectively a morphism $\phi$ as
\[
\supp A = \{ x \in X \; | \; A_x \neq 0 \} \quad \mbox{ and } \quad
\supp \phi = \{ (x,y) \in X \times X \; | \; \phi_{y,x} \neq 0 \}.
\]
If a set $\cale$ of subsets of $X \times X$ and a set $\calf$ of
subsets of $X$
satisfies the conditions (i)-(iv) listed in
\cite[Subsection~2.3]{BFJR-TOP}
then we speak of $\cale$ respectively $\calf$ as morphism and
object support conditions and define
\[
\calc ( X ; \cale , \calf ; \cala )
\]
to be the subcategory of $\calc ( X ; \cala)$ consisting of
objects $A$ for which there exists
an $F \in \calf$ such that $\supp A \subset F$ and morphisms
$\phi$ for which there exists an $E \in \cale$
such that $\supp \phi \subset E$.
Observe that for $g \in G$ we have
\[
\supp g^{\ast}A = g^{-1} (\supp A) \quad \mbox{and} \quad
\supp g^{\ast}\phi = g^{-1} (\supp \phi).
\]
We say that $\cale$ is $G$-invariant if for every $g \in G$,
$E \in \cale$ we have $g(E) \in \cale$, where $G$ acts diagonally on
$X \times X$.
We say that $\calf$ is $G$-invariant if for every $g \in G$,
$F \in \calf$ we have $g(F) \in \calf$.
For $G$-invariant object and morphism support
conditions $\cale$ and $\calf$ there is a  $G$-action on
$\calc ( X ; \cale , \calf ; \cala )$ and we can consider the
corresponding fixed category, which we denote
\[
\calc^G( X ; \cale , \calf ; \cala ).
\]
The following example was used in
Proposition~\ref{prop:delooping-of-assembly-map}.

\begin{example}
\label{ex:lambda-cala}
Let $X = [0,\infty)$.
Let $\cale = \{ E_\alpha  \; | \; \alpha > 0 \}$, where
$E_\alpha = \{ (x,y) \in [0,\infty)^{\x 2}
                         \; | \; |x-y|< \alpha \}$
and $\calf = \{ [0,r] \; | \; r \in \IR \}$.
Then there is an Eilenberg swindle on
$\Lambda \cala = \calc ( [0,\infty); \cale ; \cala )$
induced by the map $t \mapsto t+1$ on $[0,\infty)$ and
$\cala ' = \calc( [0,\infty); \cale, \calf ; \cala ) \subset
 \Lambda \cala = \calc ( [0,\infty); \cale ; \cala )$
is a Karoubi filtration, see
\cite[Definition 1.27]{Carlsson-Pedersen-Controlled-algebra-Novikov}.
\end{example}

\subsection{Assembly as Forget-Control}

From this point on it is clear that every construction and every proof in \cite{BFJR-TOP}
which treats the category of finitely generated $R$-modules in a formal way does have an
analog in our context. In particular there is a category
\[
\cald^G ( X ; \cala)
\]
defined analogously to the category $\cald^G ( X)$ from Subsection~3.2 in \cite{BFJR-TOP}
and this construction is functorial in the $G$-space $X$. The functor
\[
X \mapsto \IK^{-\infty}  \cald^G ( X ) .
\]
is a $G$-equivariant homology theory on the category of $G$-CW complexes, compare \cite[Section~4]{BFJR-TOP}.

We will now identify this controlled version of a $G$-equivariant homology theory with the $G$-equivariant homology theory
that we defined in Section~\ref{sec:FJ-with-coeff} via the $\Or G$-spectrum $\bfK_\cala$ from Definition~\ref{def:KuntenA}.

\begin{theorem}
\label{thm:identify-assembly-map}
There is an isomorphism between the functors
$X \mapsto H^{G}_*(X;\bfK_\cala)$ and
$X \mapsto \pi_{*+1}(\IK^{-\infty} \cald(X;\cala)^G)$
from $G$-$CW$-complexes to graded abelian groups.
In particular, the map
\[
K_{*+1}(\cald(E_\calf G, \cala)^G) \to K_{*+1}(\cald(\pt, \cala)^G)
\]
is a model for the assembly map
\[
H_*^{G}(E_\calf G, \bfK_\cala) \to H_*^{G}(\pt, \bfK_\cala).
\]
\end{theorem}

\begin{proof}
Without twisted coefficients this was done in
\cite[Section~6]{BFJR-TOP}.
The proof in the case with
twisted coefficients is essentially the same.
The only step in the proof
where the argument needs to be rethought is
Step (ii) in the proof of \cite[Proposition~6.2]{BFJR-TOP}.
This step is redone in lemma~\ref{lem:comparison} below.
\end{proof}

\begin{lemma} \label{lem:comparison}
Let $T$ be a $G$-set.
Let $\calf_{Gc}$ the object support condition on $T \x G$ that
contains exactly the $G$-compact subsets.
Let $\cale_\Delta$ be the morphism control condition on $T$
that contains only the diagonal of $T$.
Let $p \colon T \x G \to T$ denote the projection.
There is an additive functor
\begin{eqnarray*}
F \colon \cala \ast_G T \to
\calc ( T \times G , p^{-1} \cale_{\Delta} , \calf_{Gc} ; \cala)^G
\end{eqnarray*}
which yields an equivalence of categories.
This functor is natural in $T$.
\end{lemma}

\begin{proof}
It is straightforward to check that
\[
F(A)_{(t, g^{-1} )}  = g^{\ast} ( A_{gt} ), \quad
F( \phi )_{(t^{\prime} , k^{-1}),(t , g^{-1} )} =
\left\{
\begin{array}{cl} g^{\ast} ( \phi_{kg^{-1} ,gt} ) &
\mbox{ if } t^{\prime}=t  \\
              0
& \mbox{ if } t^{\prime} \neq t \end{array}
\right.
\]
defines an additive functor. Here,
in order to check that $F(A)$ satisfies the object support
condition observe that
$G(t,g) \mapsto g^{-1}t$ is a  bijection between the orbits of
the left $G$-set $T \times G$ and the set $T$.
Because of the object support condition an object in the
target category can be written as a direct sum of objects supported on
a single orbit of $T \times G$.
Because of the $G$-invariance an
object $C=(C_{(t,g)})$ supported on a single orbit $G(t,g)$ is determined
by its value at one point of the orbit together with
the $G$-action on the category $\cala$. Now
the object $A \in \cala \ast_G T$ supported on the single
point $\{ t\}$
which is given by $A_t=C_{(t,e)}$ maps
to $C$ under the functor $F$. Since the functor is additive
we conclude that every object in the target category is isomorphic
to an object in the image of the functor $F$.
The functor is easily seen to be faithful. It remains to prove
that it is full, i.e.\ surjective on morphism sets. If
$f = (f_{(t^{\prime},k^{-1}),(t , g^{-1} )} )$ is a morphism in
the target category then because of the
$p^{-1}\cale_{\Delta}$-condition
$f_{(t^{\prime},k^{-1})},(t , g^{-1} )$ is
non-trivial only if $t^{\prime}=t$. If one
defines a morphism $\phi$ in $\cala \ast_G T$ by
$\phi_{k,t} = f_{(t,k^{-1})(t,e)}$ then one can use the
$G$-invariance of $f$ in order to check that $F( \phi )= f$.
\end{proof}

\section{Applications}
\label{sec:Applications}

As already mentioned many arguments in controlled algebra treat the
category $\cala$ as a formal variable. Consequently existing proofs
for results about the Farrell-Jones Conjecture without coefficients can
often be carried over to the context with coefficients.
We will state three results obtained in this way.

The following is a generalization of the main
Theorem in \cite{BR-JAMS}.

\begin{theorem}
\label{thm:BR-with-coefficients}
Let $G$ be the fundamental group of a closed Riemannian manifold
of strictly negative sectional curvature.
Then the algebraic $K$-theory Farrell-Jones Conjecture
with Coefficients~\ref{con:FJC-with-coefficients-alg-K-theory}
holds for $G$.
\end{theorem}
\begin{proof}
Even though the original proof (injectivity in \cite{BFJR-TOP} and surjectivity in \cite{BR-JAMS})
is quite lengthy, one can quickly check that
the only places in the proof where the arguments need to be rethought are the following.
\begin{enumerate}
\item
The functor $\ind$ appearing in Proposition~8.3 in \cite{BFJR-TOP} needs to be
promoted to a functor $\ind \colon \cald^H( X ; \res_H \cala ) \to \cald^G ( G \times_H X ; \cala )$,
where $X$ is an $H$-space, $H$ is a subgroup of $G$ and $\cala$ is an additive category with $G$-action.
The new formulas for the functor $\ind$ are
\[
(\ind M)_{[g,x]} = ( g^{-1} )^{\ast} M_x \quad \mbox{ and }
      \quad (\ind \phi)_{[g^{\prime},x],[g,x]} = (g^{-1})^{\ast} ( \phi_{ g^{-1}g^{\prime} x^{\prime} , x } )
\]
if $g^{-1}g^{\prime} \in H$ and $0$ otherwise.
\item
The proof of injectivity uses the injectivity result for the assembly with
respect to the trivial family from \cite{Carlsson-Pedersen-Controlled-algebra-Novikov}, compare
(iii) in Subsection~10.3 in \cite{BFJR-TOP}.
We hence need the version with coefficients of that result.
It is a special case of Theorem~\ref{thm:david-with-coefficients} below.
\item
In order to define the ``Nil''-spectra denoted $\IN_i$ in Subsection~10.2 in \cite{BFJR-TOP}
one needs that the assembly map for the infinite
cyclic group with respect to the trivial family is split injective.
This is the special case of Theorem~\ref{thm:david-with-coefficients} below, where $G$ is the infinite cyclic group.
\item
The construction of the transfer functor and the proof of its properties in Section~5 of \cite{BR-JAMS}
need to be adapted to the set-up with coefficients.
This will occupy the rest of this proof.
\end{enumerate}
It will be convenient to restrict this discussion to connective
$K$-theory because we use Waldhausen categories.
This suffices by Corollary~\ref{cor:inheritance-n-to-n-1}.
(On the other hand, this discussion can be extended to
non-connective $K$-theory, by giving an adhoc definition of
non-connective $K$-theory for the Waldhausen categories we
encounter in the following,
compare \cite[Remark 5.3]{BR-JAMS}.)

First we need a replacement for the category of homotopy finite chain complexes, defined in Subsection~5.2 and 8.1 in \cite{BR-JAMS}.
Given a category $\cala$ with a right $G$-action and an infinite cardinal number $\kappa$ (chosen large enough)
we construct below in Lemma~\ref{lem-constr-A-kappa} a category with right $G$-action
$\cala^{\kappa}$. Analogously to $\overline{\calc}^G (X ; \cale )$ from \cite{BR-JAMS}
we define $\overline{\calc}^G ( X ; \cale ; \cala^{\kappa} )$ by allowing objects $M =  (M_x )_{x \in X}$, where
$M_x$ is an object in $\cala^{\kappa}$ and the support of $M$ is an arbitrary subset of $X$.
The category $\overline{\calc}^G ( X ; \cale ; \cala^{\kappa} )$
plays the role of the category that is (unfortunately) called $\overline{\cala}$ in
Subsection~8.1 in \cite{BR-JAMS}. Hence the category
$\ch_{hf} \calc^G ( X ; \cale ; \cala )$ is defined to be the ``homotopy closure'' of
the category $\ch_f \calc( X ; \cale ; \cala)$ inside
$\ch \overline{\calc}^G( X ; \cale ; \cala^{\kappa} )$.
The fibre complex $F$ and its variants from Subsection~5.3 in  \cite{BR-JAMS} can be considered as objects
in $\ch \overline{\calc}^G ( \tilde{E} \times \IT ; \cale ; \calf^{\kappa}(\IZ))$, which is defined analogous to
$\ch \overline{\calc}^G ( X ; \cale ; \cala^{\kappa} )$.
Here $\calf^{\kappa}(\IZ)$ denotes a small model for the category of those free $\IZ$-modules which admit
a basis of cardinality less than or equal to $\kappa$. The category $\calf^{\kappa}(\IZ)$ carries the trivial $G$-action.
The ``tensor product'' $- \otimes - \colon \cala^{\kappa} \times \calf^{\kappa} ( \IZ ) \to \cala^{\kappa}$
from Lemma~\ref{lem-constr-A-kappa}~\ref{lem-sub-tensor}
now allows to construct the
the transfer functor $M \mapsto M \otimes F$, $\phi \mapsto \phi \otimes \nabla$ as before.
The proof carries over without change until the end of Subsection~5.4 in \cite{BR-JAMS}.
In Proposition~5.9 the action of the Swan group $\Sw ( G ; \IZ )$ on $K_n ( RG )$ needs
to be replaced with the action of the Swan group on
$K_n ( \cala \ast_G \pt )$ that we describe below in Section~\ref{sec-swan-group-action}.
The proof of Proposition~5.9 in \cite{BR-JAMS} remains unchanged
until one reaches diagram~(5.12). That diagram now gets replaced with the following diagram
\[
\xymatrix{
\calc^G ( Gb_0 ; \cala)
\ar[r]^-{- \otimes \nabla}  &
\ch_{hf} \calc^G ( Gb_0 ; \cala ) \\
\cala \ast_G \pt \ar[u]^-{F}
\ar[r]^-{- \otimes F_0}  &
\ch_{hf} (\cala \ast_G \pt) \ar[u]^-{\ch_{\hf} F}.
         }
\]
Here $\ch_{\hf} ( \cala \ast_G \pt)$ is the homotopy closure of $\ch_f (\cala \ast_G \pt)$
in $\ch (\cala^{\kappa} \ast_G \pt)$.
The functor  $F$ is a special case of the equivalence from Lemma~\ref{lem:comparison}
and is given by $F(A)_{g^{-1}b_0} = g^{\ast} (A)$ and $F( \phi)_{k^{-1}b_0 , g^{-1}b_0} = g^{\ast} ( \phi_{kg^{-1}} )$.
The functor $- \otimes \nabla$ is given by $(M_{gb_0} ) \otimes \nabla = ( M_{gb_0} \otimes F_0 )$
and $( \phi_{gb_0 , hb_0} ) \otimes \nabla = ( \phi_{gb_0 , hb_0} \otimes \nabla_{gb_0 , hb_0 } )$.
We can equip $F_0$ with a $G$-action in such way that $\nabla_{gb_0 , hb_0} \colon F_0 \to F_0$
corresponds to $l_{gh^{-1}}$, i.e.\ to left multiplication with $gh^{-1}$.
With this notation
the functor $- \otimes F_0$ is defined as follows.
The object $A$ maps to $A \otimes F_0$ the morphism $\phi = ( \phi_g )$ maps to
$(\phi_g \otimes l_g)$.
As opposed to the original diagram in \cite{BR-JAMS} the diagram now commutes.
Let $\inc \colon \cala \ast_G \pt \to \ch_{hf} (\cala \ast_G \pt)$ denote the inclusion.
It follows from the discussion of the Swan group action in Section~\ref{sec-swan-group-action}
below that on the level of $K$-theory $\inc^{-1} \circ (- \otimes F_0)$ corresponds to
multiplication with $[F_0]= \sum_i (-1)^i [H_i (F_0)] \in \Sw^{\ch}( G ; \IZ ) \cong \Sw( G ; \IZ )$.
\end{proof}

The following is a generalization of
a result of Rosenthal \cite{Rosenthal-splitting-K-theory}.

\begin{theorem}
\label{thm:david-with-coefficients}
Let $G$ be a group.
Suppose that there is a model for $E_\Fin G$
that is a finite $G$-$CW$-complex and admits
a compactification $X$ such that
\begin{enumerate}
\item the $G$-action extends to $X$;
\item $X$ is metrizable;
\item $X^F$ is contractible for every $F \in \Fin$;
\item $(E_\Fin G)^F$ is dense in $X$ for every $F \in \Fin$;
\item compact subsets of $E_\Fin G$ become small near
      $Y = X - E_\Fin G$.
      That is, for every compact subset of $E_\Fin G$
      and for every neighborhood $U \subset X$ of $y \in Y$,
      there exists a neighborhood $V \subset X$ of $y$ such that
      $g \in G$ and $gK \cap V \neq \emptyset$ implies
      $gK \subset U$.
\end{enumerate}
Let $\cala$ be an additive
category with right $G$-action.
Then the assembly map
\[
H^G_*(E_\Fin G; \bfK_\cala) \to H^G_*(\pt;\bfK_\cala)
\]
is split injective.
\end{theorem}

The following is a generalization of
the main result from \cite{Bartels-finite-asymp-dim}.

\begin{theorem}
\label{thm:bartels-fin-asy-with-coeff}
Let $G$ be group of finite asymptotic dimension
that admits a finite model for the classifying space $BG$.
Let $\cala$ be an additive category with right
$G$-action.
Then the assembly map
\[
H^G_*(EG;\bfK_\cala) \to H^G_*(\pt;\bfK_\cala)
\]
is split injective.
\end{theorem}

Both, the proof of Theorem~\ref{thm:david-with-coefficients} and
of Theorem~\ref{thm:bartels-fin-asy-with-coeff} are trivial modifications of the original proofs
in \cite{Rosenthal-splitting-K-theory} respectively
\cite{Bartels-finite-asymp-dim}.
Everywhere in these proofs the category of $R$-modules is treated as a formal variable and
can simply be replaced by the additive category with $G$-action $\cala$.

\section{The Swan group action} \label{sec-swan-group-action}
\label{sec:swan-group-action}

In the proof of Theorem~\ref{thm:BR-with-coefficients} we used some facts about Swan group actions which we are going to prove now.
The Swan group $\Sw( G ; \IZ)$ is the $K_0$-group of the category of $\IZ G$-modules that are finitely generated as $\IZ$-modules.
Recall from Subsection~8.2 of \cite{BR-JAMS} that there are version $Sw^{\fr}(G ; \IZ)$ and $\Sw^{\ch}( G; \IZ )$ defined using
$\IZ G$-modules that are finitely generated free as $\IZ$-modules,
respectively bounded below chain complexes of $\IZ G$-modules that are free as $\IZ$-modules
and whose homology is finitely generated as a $\IZ G$-module. We have shown in Proposition~8.3 in \cite{BR-JAMS} that the natural maps
$j \colon \Sw^{\fr} ( G ; \IZ ) \to \Sw^{\ch} ( G ; \IZ )$ and $i \colon \Sw^{\fr} ( G ; \IZ ) \to \Sw ( G ; \IZ )$  are isomorphisms.

\begin{lemma}
Let $\cala$ be an additive category with right $G$-action.
There exists for $n \geq 1$ a commutative diagram
\[
\xymatrix{
  K_n ( \cala \ast_G \pt ) \otimes_{\IZ}
    \Sw^{\fr} (G ; \IZ )  \ar[d]^-{\id \otimes j } \ar[r]
  &
  K_n ( \cala \ast_G \pt ) \ar[d]^-{\inc}
  \\
  K_n ( \cala \ast_G \pt ) \otimes_{\IZ}
     \Sw^{\ch} (G ; \IZ ) \ar[r]
  &
  K_n ( \ch_{\hf} (\cala \ast_G \pt) ) ,
}
\]
where both vertical arrows are isomorphisms.
In this way $K_n ( \cala \ast_G \pt )$ becomes a module
over the Swan ring
$\Sw^{\fr}( G ; \IZ )\cong \Sw^{\ch} ( G ; \IZ)
                             \cong \Sw (G ; \IZ )$.
\end{lemma}

Using the suspension category $\Sigma \cala$ from
Proposition~\ref{prop:delooping-of-assembly-map} it
is possible to formulate a version of the above Lemma
that applies to all $n$.
However, for our purposes the above formulation suffices.

\begin{proof}
Choose a small model $\calf^f(\IZ)$ for the category of finitely generated free $\IZ$-modules such that the underlying $\IZ$-module of
every $\IZ G$-module that is used in the construction of $\Sw^{\fr} (G ; \IZ)$ is contained in $\calf^f ( \IZ)$.

We replace $\cala$ by the equivalent category  $\cala^f$ from Lemma~\ref{lem-constr-A-kappa} but refer to it as $\cala$ in the following.
This is justified because $- \ast_G \pt$ respects equivalences, compare Remark~\ref{rem:cala-ast-T-is-natural}.
We can hence assume
that there exists a tensor product $- \otimes - \colon \cala \times \calf^f( \IZ ) \to \cala$ with the expected properties.
Then define for a Swan module $M$ the additive functor $- \otimes M \colon \cala \ast_G \pt \to \cala \ast_G \pt$ by
\[
A \otimes M = A \otimes UM, \quad \mbox{ and } \quad (\phi \otimes M )_{g } = \phi_g \otimes l_g .
\]
Here $UM$ denotes the underlying $\IZ$-module of $M$ and $l_g$ denotes left multiplication by $g$.
Note that the $\IZ G$-module structure of $M$ enters only in the morphisms.
A morphism $f \colon M \to N$ of Swan modules induces a natural transformation $\tau(f,A) \colon A \otimes M \to A \otimes N$ that
is given by $\tau(f,A)_g = 0$ if $g \neq e$ and $\tau(f,A)_e = \id_A \otimes f$. A short exact sequence
$L \to M \to N$ of Swan-modules leads to
a short exact sequence of functors, because the underlying sequence $UL \to UM \to UN$ always splits and being a short exact sequence
of functors is checked objectwise.
One uses
\cite[1.3.2~(4)]{WaldhausenLNM1126} in order to check that one obtains the $\Sw^{\fr} ( G ; \IZ )$-action.

For the $\Sw^{\ch} ( G ; \IZ )$-action arrange that the underlying $\IZ$-modules of all
Swan modules that appear in the chain complexes lie in a small model
$\calf^{\kappa} (\IZ)$ of the category of finitely generated free modules that
admit a basis $B$ of cardinality $\card(B) \leq \kappa$.
(Strictly speaking we have to have a cardinality assumption when we define the category of chain
complexes that leads to $\Sw^{\ch} ( G ; \IZ )$.)
By Lemma~\ref{lem-constr-A-kappa} there exists an inclusion $\cala \to \cala^{\kappa}$ and the
``tensor product'' we used so far extends to a tensor product
$- \otimes - \colon \cala^{\kappa} \otimes \calf^{\kappa} ( \IZ) \to \cala^{\kappa}$.

Define $\ch_{\hf} (\cala \ast_G \pt)$ to be the category of
chain complexes in $\ch (\cala^{\kappa} \ast_G \pt)$
that are chain homotopy equivalent to a bounded below and above
chain complex in $\ch (\cala \ast_G \pt)$.
Similarly let $\ch_{\hf} \calf^f( \IZ ) $ denote the
category of chain complexes in $\ch \calf^{\kappa} ( \IZ)$
which are chain homotopy equivalent
to a finite complex in $\ch \calf^f( \IZ )$,
compare \cite[Section~8.1]{BR-JAMS}.
The right hand vertical arrow induced by the inclusion
$\cala \ast_G \pt \to \ch_{\hf} (\cala \ast_G \pt)$
is an equivalence by \cite[Lemma~8.1]{BR-JAMS}.
Now for a chain complex $C_{\bullet}$
which represents an element in $\Sw^{\ch} (  G ; \IZ)$
one defines a functor
\[
- \otimes C_{\bullet} \colon \cala \ast_G \pt
                   \to \ch_{\hf} (\cala \ast_G \pt)
\]
analogously to $- \otimes M$ above by $A \otimes C_{\bullet} = A \otimes UC_{\bullet}$.
This is well defined because $UC_{\bullet}$ is by definition of $\Sw^{\ch} ( G ; \IZ)$
a bounded below complex of finitely generated free $\IZ$-modules whose homology is concentrated
in finitely many degrees and each of its homology groups is a finitely generated $\IZ$-module.
Such a complex is homotopy equivalent to a complex $C^{\prime}_{\bullet}$
of finitely generated free $\IZ$-modules which is concentrated in finitely many degrees.
(In order to prove this assume that $C_{\bullet}$ is concentrated in non-negative degrees and
use induction over the largest number $m$ such that $H_m( C_{\bullet} ) \neq 0$.
In the case $m=0$ the complex $UC_{\bullet}$ is a resolution of the finitely generated $\IZ$-module
$H_0 ( C_{\bullet} )$ and is hence homotopy equivalent to a finite resolution.
For $m \geq 1$ choose a finite resolution $D_{\bullet}$ of $H_m ( C_{\bullet} )$ and
construct a $H_m$-isomorphism $f \colon D_{\bullet} \to UC_{\bullet}$.
Factorize $f$ over its mapping cylinder $\cyl(f) \simeq UC_{\bullet}$ and study the
sequence $D_{\bullet} \to \cyl(f) \to \cone(f)$.)
Consequently $A \otimes UC_{\bullet}$ is homotopy equivalent to $A \otimes C^{\prime}_{\bullet}$
and hence lies in $\ch_{\hf} \cala \ast_G \pt$.
Again a short exact sequence of chain complexes leads to a short exact sequence of functors, because the objects
depend only on the underlying $\IZ$-chain complexes.
A homology equivalence
$C_{\bullet} \to D_{\bullet}$ can be considered as a homotopy equivalence $UC_{\bullet} \to UD_{\bullet}$
and hence induces a homotopy equivalence
$A \otimes UC_{\bullet} \to A \otimes UD_{\bullet}$.
\end{proof}
In the proof above and in the proof of Theorem~\ref{thm:BR-with-coefficients} we used the following lemma.

\begin{lemma} \label{lem-constr-A-kappa}
Let $\cala$ be a small additive category with a $G$-action
by additive functors.
Let $\kappa$ be a fixed infinite cardinal.
Denote by $\calf^{\kappa} ( \IZ )$ some small model for the category of all free $\IZ$-modules
which admit a basis $B$ with $\card(B)\leq \kappa$.
Equip $\calf^{\kappa}(\IZ)$ with a tensorproduct functor $- \otimes_{\IZ} -$.
(This of course involves choices.)
Let $\calf^f(\IZ)$ be the full subcategory of $\calf^{\kappa}(\IZ)$ that consists of finitely generated free $\IZ$-modules.

There exist additive categories $\cala^f$ and  $\cala^{\kappa}$
with $G$-action and
$G$-equivariant additive inclusion functors
\[
\cala \to \cala^f \to \cala^{\kappa}
\]
such that the following conditions hold.
\begin{enumerate}
\item
The inclusion $\cala \to \cala^f$ is an equivalence of categories.
\item
In $\cala^{\kappa}$ there exist categorical sums over indexing sets $J$ with $\card(J) \leq \kappa$. Hence if  $A_j$, $j \in J$ is
a family of objects in $\cala^{\kappa}$ then $\bigoplus_{j \in J} A_j$ exists.
Moreover one can make a choice for these sum objects such that
for every $g \in G$ we have an equality
\[
g^{\ast} ( \bigoplus_{j \in J} A_j ) = \bigoplus_{j \in J} g^{\ast} ( A_j ).
\]
\item \label{lem-sub-tensor}
There exists a bilinear bifunctor
\[
- \otimes - \colon \cala^{\kappa} \times \calf^{\kappa} ( \IZ ) \to \cala^{\kappa},
\]
which restricts to
\[
- \otimes - \colon \cala^f \times \calf^f ( \IZ ) \to \cala^f.
\]
The functor
is compatible with direct sums in the sense that for a
family of objects $A_j$, $j \in J$ in $\cala^{\kappa}$ with $\card(J) \leq \kappa$
and a $\IZ$-module $F \in \calf^{\kappa} ( \IZ ) $ there exists a natural isomorphism
\[
\left( \bigoplus_{j \in J } A_j \right) \otimes F \cong  \bigoplus_{j \in J } \left(A_j  \otimes F \right)
\]
(The analogous statement for finite direct sums holds in both variables because of the bilinearity.)
The bifunctor is also compatible with the $G$-action in the sense that for every $g \in G$
and all morphisms $\phi \colon A \to B$ in $\cala^{\kappa}$
and $f \colon F \to F^{\prime}$ in $\calf^{\kappa} ( \IZ ) $ we have equalities
\[
g^{\ast} ( A \otimes F ) = ( g^{\ast}A ) \otimes F \quad \mbox{ and } \quad g^{\ast} ( f \otimes \phi ) = g^{\ast} ( f ) \otimes \phi.
\]
\item
For objects $A$ in $\cala^{\kappa}$ and $F$, $F^{\prime}$ in $\calf^{\kappa} ( \IZ ) $ we have a natural isomorphism
\[
(A \otimes F) \otimes F^{\prime} \cong A \otimes (F \otimes F^{\prime}) .
\]
\end{enumerate}
\end{lemma}

\begin{proof}
The construction of the categories $\cala^f$ and $\cala^{\kappa}$ makes use of the following two elementary constructions.
First construction: let $\calb$ and $\calc$ be two $\Ab$-categories,
i.e.\ categories enriched over abelian groups, compare \cite[I.8]{MacLane-cat-for-working}.
Then we define the $\Ab$-category $\calb \otimes \calc$ as follows.
Objects are pairs of objects which we denote $B \otimes C$, where $B$ is an object in $\calb$ and $C$ an object in $\calc$.
We set
\[
\mor_{\calb \otimes \calc} ( B \otimes C  , B^{\prime} \otimes C^{\prime} )
        = \mor_{\calb} (B, B^{\prime} ) \otimes_{\IZ} \mor_{\calc} (C , C^{\prime} ).
\]
Composition and identities are defined in the obvious way.
The construction is functorial with respect to additive
functors in  $\calb$ and $\calc$ and preserves
additive equivalences.

Second construction: given an $\Ab$-category $\cald$ and a set $I$ we
define the category $\cald(I)$ to be the category whose objects are
families $D=(D(i))_{i \in I}$ of objects in $\cald$ and where a morphism $f \colon D \to D^{\prime}$ is a family of morphisms
$f(j,i) \colon D(i) \to D(j)$, $i$, $j \in I$ subject to the condition that for a fixed $i \in I$ there are only finitely many
$j \in I$ such that $f(j,i) \neq 0$. Composition is the usual matrix multiplication,
where the components of $f^{\prime} \circ f$ are given by
\[
(f^{\prime} \circ f )(k,i) = \sum_j f^{\prime} ( k,j) \circ f(j,i).
\]

Now set $\cala^f = \cala \otimes \calf^f( \IZ)$. The inclusion $\cala \to \cala^f$ is given
by $A \mapsto A \otimes \IZ^1$, where $\IZ^1$ is some $1$-dimensional free $\IZ$-module in $\calf^f( \IZ )$.
It is not difficult to check that this is an equivalence of categories.
We define a  ``tensor product'' $- \otimes - \colon \cala^f \times \calf^f ( \IZ ) \to \cala^f$
by $(A \otimes F) \otimes F^{\prime} = A \otimes ( F \otimes F^{\prime} )$.
The $G$-action on $\cala^f$ is defined by $g^{\ast} ( A \otimes F ) = (g^{\ast}A ) \otimes F$.

Next we choose a set $I$ of cardinality $\kappa$ and set
\[
\cala^{\kappa} = \left( \cala^f \otimes \calf^{\kappa} ( \IZ) \right) (I).
\]
The inclusion functor $\cala^f \to \cala^{\kappa}$ sends $A \otimes F$ to the object
which at some fixed index $i_0$ is given by $A \otimes F$ and is zero everywhere else.
The $G$-action extends via $g^{\ast} ( (A(i) \otimes F(i))_{i \in I} ) = ( (g^{\ast}A(i)) \otimes F(i) )_{i \in I}$.
The ``tensor product'' extends by
$( (A(i) \otimes F(i))_{i \in I} ) \otimes F =  (A(i) \otimes (F(i) \otimes F)  )_{i \in I}  $, where of course
the $A(i)$ are objects in $\cala$ and the $F(i)$ and $F$ are objects in $\calf^{\kappa}(\IZ)$.
The existence of the required direct sums is a consequence of the fact
that $\card(I \times I) = \card (I)$ for an infinite set $I$, see for example \cite[Appendix~2, \S~3, Theorem~3.6]{Lang}.
\end{proof}


\end{document}